\documentclass[12pt]{article}
\usepackage{amsfonts}
\usepackage[pctex32]{graphics}
\usepackage[dvips]{graphicx}
\usepackage{amsthm}
\usepackage{amsmath}
\usepackage{amssymb}
\usepackage{fontenc}
\usepackage[all]{xy}
\usepackage{latexsym}
\usepackage{layout}
\usepackage{epsfig}
\usepackage{graphicx}
\usepackage{pstricks,pst-node,pst-plot} 
\usepackage{color} 
\title{Telescoping Sums, Permutations, and First Occurrence Distributions}
\author {Anant Godbole and Jie Hao\\ East Tennessee State University}
\begin{document}
\def\qed{\vbox{\hrule\hbox{\vrule\kern3pt\vbox{\kern6pt}\kern3pt\vrule}\hrule}}
\def\ms{\medskip}
\def\n{\noindent}
\def\ep{\varepsilon}
\def\G{\Gamma}
\def\lr{\left(}
\def\ls{\left[}
\def\rs{\right]}
\def\lf{\lfloor}
\def\rf{\rfloor}
\def\lg{{\rm lg}}
\def\lc{\left\{}
\def\rc{\right\}}
\def\rr{\right)}
\def\ph{\varphi}
\def\p{\mathbb P}
\def\nk{n \choose k}
\def\a{\cal A}
\def\s{\cal S}
\def\e{\mathbb E}
\def\l{\lambda}
\def\v{\mathbb V}
\newcommand{\bsno}{\bigskip\noindent}
\newcommand{\msno}{\medskip\noindent}
\newcommand{\oM}{M}
\newcommand{\omni}{\omega(k,a)}
\newtheorem{thm}{Theorem}[section]
\newtheorem{con}{Conjecture}[section]
\newtheorem{claim}[thm]{Claim}
\newtheorem{deff}[thm]{Definition}
\newtheorem{lem}[thm]{Lemma}
\newtheorem{cor}[thm]{Corollary}
\newtheorem{rem}[thm]{Remark}
\newtheorem{prp}[thm]{Proposition}
\newtheorem{ex}[thm]{Example}
\newtheorem{eq}[thm]{equation}
\newtheorem{que}{Problem}[section]
\newtheorem{ques}[thm]{Question}
\providecommand{\floor}[1]{\left\lfloor#1\right\rfloor}
\maketitle
\date{}
\section{Telescoping Series}  Telescoping sums are a delight.  While Euler and many others since have  produced magnificent proofs of the fact that
\[\zeta(2)=\sum_{n=1}^\infty\frac{1}{n^2}=\frac{\pi^2}{6},\] 
the ``related" telescoping sum
\[\sum_{n=1}^\infty\frac{1}{n(n+1)}=\lr1-\frac{1}{2}\rr+\lr\frac{1}{2}-\frac{1}{3}\rr+\ldots\]
can be easily seen to telescope to 1.  Telescoping sums are also one of the few classes of infinite series for which one can go beyond the somewhat unsatisfactory Calculus 2 assertion that ``this series converges," by exhibiting convergence to a specific real number.  Questions involving telescoping sums are often to be found on the Putnam and other Exams \cite {ga}, where algebraic and trigonometric identities work in tandem with the underlying telescoping nature of the sum to produce stellar formulas for finite and infinite sums.  Telescoping series are those, in fact, for which the sum is revealed through an ``antiderivative" method akin to the one used when employing the Fundamental Theorem of Calculus.  In this note, we study the simplest type of infinite telescoping sum, namely
\begin{equation}\sum_{n=1}^\infty p_n=\sum_{n=1}^\infty (q_n-q_{n+1})=1,\end{equation}
due to our assumptions that $1=q_1\ge q_2\ge q_3\ldots\ge0$ and $\lim_{n\to\infty}q_n=0$.  Now (1) clearly defines a probability distribution on ${\mathbb Z}^+$, but making unmotivated choices such as $q_n=\frac{1}{3n^4-2n^2}$ would serve little purpose, however, so let us state our objective:  We wish to propose three probability distributions that each 

(i)  are similar in form to familiar discrete distributions; 

(ii) arise from the cycle (\cite{b}) or pattern containment (\cite{b}, \cite{k}) structure or random permutations; and 

(iii) lead to open-ended sets of questions.  

\noindent Let us go one step further, by asking 

(iv) that there be events $A_n,n\ge1$ in some probability space so that $q_n=\p(A_n)$ with $A_{n+1}\subseteq A_n$ so that, as we shall see, $p_n,n\ge 1$ provides a model for {\it first occurrence distributions} of certain events.  

The reader should note that discrete models, other than those provided by random permutations, could equally well have been used as launching pads for our discrete distributions, and s/he is invited to  come up with well-motivated first occurrence distributions along the lines of the ones in this paper.

\section{Discrete Distributions on ${\mathbb Z}^+$}  Our motivating distributions on ${\mathbb Z}^+$ are well-known ones that are studied in elementary probability texts such as \cite{r}:  The Zeta, the Poisson, and the Geometric.

The zeta distribution with parameter $k>1$ has associated probability mass function (pmf) given by 
\[f(n)=\frac{1}{\zeta(k)}\cdot\frac{1}{n^k}, n=1, 2, \ldots,\]
where $\zeta(k)=\sum_{n=1}^\infty\frac{1}{n^k}$ is the Riemann zeta function.  Early uses of the zeta distribution were by Zipf and Pareto, the latter in the context of income distributions \cite{r}.  For this reason the zeta distribution is also known as the Zipf or discrete Pareto distribution.  In recent years \cite{clv}, the zeta distribution has been referred to as the power law distribution, and been been used to model vertex degrees in ``small world" networks, notably when $k\in(2,3)$.

The Poisson distribution with parameter $\l>0$, denoted by Po($\l$), is the distribution of a random variable $X$ with  pmf
\[f(x)=\frac{e^{-\l}\l^x}{x!}, x=0,1,\ldots.\]
First studied in the context of tallies of Prussian soldiers' deaths per year by kicks from horses, Po($np$) is a wonderful approximation for the Binomial distribution with parameters $n$ and $p$ when $p$ is small -- earning it the moniker of a ``rare event distribution."  The second, equivalent, way that Po($\l$) arises is as the count of events in time or space that satisfy some mild conditions \cite{r}.  It is also known that if $X\sim {\rm Po}(\l)$, then $\e(X)=\v(X)=\l$, and the so-called moment generating function (m.g.f.) of $X$ is
\[M_X(t)=\e\lr e^{tX}\rr=\exp\{\l(e^t-1)\}.\]

The geometric random variable counts the number of independent Bernoulli trials with success probability $p$ that need to be conducted in order to get the first success.  We have, for $x\ge 1$,
\[f(x)=(1-p)^{x-1}\cdot p;\ \e(X)=\frac{1}{p};\ \v(X)=\frac{1-p}{p^2}.\]

\section{Patterns and Cycles in Permutations}  The theory of pattern avoidance in permutations is now well-established and thriving, and a survey of the many results in that area may be found in the text by Kitaev\cite{k}.  One of the earliest and most fundamental results in the field is that the number of permutations of $[n]:=\{1,2,\ldots,n\}$ in which the longest increasing sequence is of length $\le2$, the so-called 123-avoiding permutations, is given by the Catalan numbers $C_n$, which, for $n\ge0$, are given by
\[C_n=\frac{1}{n+1}{{2n}\choose{n}}\sim \frac{C\cdot 4^n}{n^{3/2}}.\] So why are permutations in which the longest increasing sequence is of size 2 or less called 123-avoiding permutations?  Here is the reason:  We say that a permutation $\pi$ contains a subpermutation (usually called a pattern) $\rho$ of length $k$ if there exist indices $i_1<i_2<\ldots<i_k$ such that $\pi(i_1), \pi(i_2), \ldots\pi(i_k)$ are in the same relative order as $\rho$.  For example if $\rho=312$, then $\pi=54213$ contains $\rho$ with $(i_1,i_2,i_3)=(1, 3, 5)$, since $\pi(i_1), \pi(i_2),$ and $ \pi(i_3)$ are the largest, smallest, and middle terms of $\pi(i_1), \pi(i_2), \pi(i_3)$.  If $\pi$ does not contain $\rho$, we say that it {\it avoids} $\rho$.  Classical bijective techniques have been used to show  that the each of the $ijk$-avoiding permutations with $\{i,j,k\}=\{1,2,3\}$ are equinumerous, and, in fact the origins of the theory of pattern avoidance can be traced to the result of Knuth \cite{kn}, who proved that a permutation $\pi\in S_n$ could be sorted with a single stack if and only if it avoided the pattern 231.  The theory of pattern avoidance is now extraordinarily rich, and there is an Annual International Conference devoted to the latest research in the field; see, e.g., {\tt www.etsu.edu/cas/math/pp2014} for the webpage of the latest, twelfth, conference.

In contrast to permutation patterns, the study of the cycle decomposition of a permutation $\pi$ has a longer history; see \cite{b}.  In particular, we recall that any permutation can be decomposed uniquely into cycles, and the elementary Cauchy formula \cite{b} states that for $\sum_{k=1}^nkj_k=n$, the number of permutations on $[n]:=\{1,2,\ldots,n\}$ with $j_k$ cycles of length $k$ is given by
\[\frac{n!}{\prod_{k=1}^nk^{j_k}j_k!}.\]
In particular, there are $(n-1)!$ permutations on $[n]$ that are unicyclic.  Arratia and Tavar\'e \cite{at} contains a distinguished study of the cycle structure of random permutations.

\section{Three New Probability Distributions and Their Properties}  

\noindent  1. THE TELESCOPING ZETA(2) DISTRIBUTION: There are $(n-1)!$ permutations on $[n]$ that are unicyclic.  Now if a permutation on $[n+1]$ is unicyclic, then its {\it reduction} to $[n]$ is also unicyclic.  For example, the permutation
\[{123456}\atop{452613},\]
written in ``two line notation," can be reduced, on eliminating the 6 in both the first and second rows and mapping 4 to 3 directly, as the unicyclic 5-permutation
\[{12345}\atop{45231};\]  a general definition of reduction can now be easily formulated.  
Let $n\ge 1$ and consider $\pi\in S_{n+1}$.  On setting $A_n$ to be the event that the reduced permutation on $[n]$ is unicyclic, we have that $A_{n+1}\subseteq A_n$, and so $\p(A_n\cap A_{n+1}^C)=\p(A_n)-\p(A_{n+1})$.  It follows that
\[\frac{(n-1)!}{n!}-\frac{n!}{(n+1)!}=\frac{1}{n(n+1)}\]
is the probability that $\pi\in S_{n+1}$ has a unicyclic reduction but we have  $\pi(n+1)=n+1$, so that $\pi$ itself is not unicyclic, and has an $n$ cycle on $[n]$ followed by a 1-cycle on $\{n+1\}$. Given a random permutation in $S_n$, we let $X=1,2,\ldots,n$ be the largest $i\le n$ such that the reduction of $\pi$ to $[i]$ is unicyclic.  We then have
\begin{equation}p_i=\p(X=i)=\frac{1}{i}-\frac{1}{i+1}=\frac{1}{i(i+1)},\ i=1,2,\ldots,n-1,\end{equation}
with \[\p(X=n)=\frac{1}{n}.\] 

\medskip

\noindent EXAMPLE:  For $n=4$ the chance $p_4$ that $\pi$ is unicyclic is $\frac{1}{4}$, as evidenced by the six permutations 2341, 2413, 3421, 3142, 4123, and 4312 respectively.  $p_3=\frac{1}{3}-\frac{1}{4}=\frac{1}{12}$, and this fact is verified by the two permutations 2314 and 3124.  The four permutations that yield $p_2=\frac{1}{6}$ are 2134, 2143, 2431, and 4132.  Finally, $p_1=\frac{1}{2}$ due to the twelve permutations 1234, 1243, 1324, 1342, 1423, 1432, 3241, 3214, 3412, 4213, 4231, and 4321.

Given an infinite sequence $X_1,X_2\ldots$ of i.i.d. uniform random variables, for some $n$, a second cycle {\it must} be introduced with probability 1, for some $n$, on the random sequence of order statistics $X_{(1)}<X_{(2)}<\ldots<X_{(n)},$ and we thus get the infinite probability distribution
\begin{equation}p_n=\frac{1}{n}-\frac{1}{n+1}=\frac{1}{n(n+1)}; n\ge 1.\end{equation}  This distribution has infinite mean, however, with $\e(X)\sim \log n$ for a permutation on $[n]$.

A question that we asked ourselves was whether the telescoping zeta distribution could arise via the theory of permutation patterns.  We came up with this scenario:  Returning to Equation (2) (or Equation (3)), we note that the probability distribution there can also be expressed, given a random permutation on $[n]$ and $i\le  n-1$, as
\[\rho_i=p_i=\frac{{n\choose {i+1}}(i-1)!(n-i-1)!}{n!}, 1\le i\le n-1; \rho_n=p_n=\frac{1}{n}\]
whose numerator we can interpret as follows:  (i) Choose any of the $(i-1)!$ unicyclic permutations on $[i]$; (ii) realize them using any of $i+1$ numbers, with the largest occupying the $(i+1)$st spot and the rest in the same relative order as that of the unicyclic permutation; and (iii) Arrange the rest of the numbers in any way possible.  We then have the first $i$ numbers being in an order-isomorphic cyclic form, the first $i+1$ numbers {\it not} yielding a cycle, and the last $n-i-1$ numbers being arbitrary.  But {\it $\rho_i$ does not yield the probability distribution of a random variable, since the various sample outcomes are not disjoint!}
Using $n=4$ as an example again, we find that the numerators of the $\rho$ ratios arise as in Table 1:
\begin{table}[h]
\caption{\label{Table 1.} Can Permutation Patterns give the telescoping zeta distribution?}
\begin{center}
\begin{tabular}{|c|c|c|c|c|c|c|}\hline\hline
$i$ & $\rho_i$ &Sample Points\\ \hline\hline
1 & $\frac{1}{2}$ & 1234, 1243, 1324, 1342, 1432, 1423, 2314, 2341, 2413, 2431, 3412, 3421 \\ \hline
2 & $\frac{1}{6}$ & 2134, 2143, 3142, 3241\\ \hline
3 & $\frac{1}{12}$ & {\bf 2314}, 3124\\ \hline
4 & $\frac{1}{4}$ & {\bf 2341, 2413, 3421, 3142}, 4123, 4312\\ \hline
None & {}& 3214, 4132, 4213, 4231, 4321\\ \hline
\end{tabular}
\end{center}
\end{table}
 We see that even though the sum of the $\rho_i$s is 1, this occurs due to double counting of five sample points and exclusion of five points.  The reader is invited to address this issue, for general values of $n$, as s/he sees fit!

\noindent 2. THE TELESCOPING POISSON DISTRIBUTION:   
The probability that a random permutation on $[n]$ has its first ascent at positions $k, k+1$ is given, for $1\le k\le n-1$, by $\frac{k}{(k+1)!}$.  To see this, choose any one of the $k+1$ elements in positions 1 through $k+1$, except for the smallest, to occupy the $k+1$st position, and then arrange the other elements in a monotone decreasing fashion.  But, to reformulate this argument via telescoping series, we let $A_k$ be the event that the permutation restricted to the first $k$ integers is monotone decreasing, note that $A_{k+1}\subseteq A_k$ and that the probability $f(k)$ that the first ascent is at $k, k+1$ satisfies, for $k\le n-1$,
\[f(k)=\p(A_k\cap A_{k+1}^C)=\frac{1}{k!}-\frac{1}{(k+1)!}=\frac{k}{(k+1)!},\] as before.  The chance that the first ascent is at position $n$ is, of course, $\frac{1}{n!}$.  Letting $n\to\infty$, we get the discrete distribution on ${\mathbb Z}^+$ given by
\begin{equation}
f(x)=\frac{x}{(x+1)!}, x=1,2,\ldots, 
\end{equation}
which is quite similar to the unit Poisson distribution on 0,1,... with mass function $f(x)=e^{-1}/x!$.  
If $X$ is the corresponding first ascent random variable generated, e.g., by a sequence $\{X_n\}_{n=1}^\infty$ of i.i.d.~uniform random variables, then
\begin{eqnarray*}\e(X)&=&\sum_{x=1}^{\infty}x\times \frac{x}{(x+1)!}\\
&=&\sum_{x=1}^{\infty}\frac{1}{(x-1)!}-\sum_{x=1}^{\infty}\frac{1}{x!}+\sum_{x=1}^{\infty}\frac{1}{(x+1)!}\\
&=&e-(e-1)+(e-1-1)\\
&=&e-1,\end{eqnarray*}
and
\begin{eqnarray*}
\v(X)&=&\e(X^2)-(\e(X))^2\\
&=&(e+1)-(e-1)^2\\
&=&e(3-e), 
\end{eqnarray*} since
\begin{eqnarray*}\e(X^2)&=&\sum_{x=1}^{\infty}x^2\times \frac{x}{(x+1)!}\\
&=&\sum_{x=2}^{\infty}\frac{1}{(x-2)!}+\sum_{x=1}^{\infty}\frac{x}{(x+1)!}\\
&=&e+1. \end{eqnarray*}
More generally, the m.g.f. of $X$ is given by
\begin{eqnarray*}\e(e^{tX})&=&\sum_{x=1}^{\infty}e^{tx}\times \frac{x}{(x+1)!}\\
&=&\sum_{x=1}^{\infty}\frac{e^{tx}}{x!}-\sum_{x=1}^{\infty}\frac{e^{tx}}{(x+1)!}\\
&=&\sum_{x=1}^{\infty}\frac{e^{tx}}{x!}-\frac{1}{e^t}\sum_{x=1}^{\infty}\frac{{(e^t)}^{x+1}}{(x+1)!}\\
&=&(e^{e^t}-1)-\frac{1}{e^t}(e^{e^t}-1-e^t)\\
&=&e^{-t}(1-e^{e^t}+e^{t+e^t}),
\end{eqnarray*}
compared to the m.g.f. $\exp\{e^t-1\}$ of the unit Poisson r.v.
A one parameter telescoping Poisson model may be defined through the formula
\[f(x)=\frac{\theta^x}{x!}-\frac{\theta^{x+1}}{(x+1)!}, x=0,1,2,\ldots; 0<\theta<1,\]
and justified through the use of monotonicity conditions for size-biased permutations.  The baseline case corresponds to $\theta=1$.  It is routine to verify, using inequalities such as $1-x\le e^{-x}$, that the above mass function assigns lower weight to $x=0$ than does the ${\rm Po}(\theta)$ variable and higher weight for any $x\ge2$.  Moreover, for $x=1$ the Poisson variable has greater mass if and only if $\theta>\ln 2$.  Furthermore, we have for the $\theta$-telescoping Poisson variable $X$,
\begin{eqnarray*}M_X(t)=\e(e^{tX})&=&\sum_{x=0}^{\infty}e^{tx}\times(\frac{\theta^x}{x!}-\frac{\theta^{x+1}}{(x+1)!})\\
&=&\sum_{x=0}^{\infty}\frac{(\theta e^t)^x}{x!}-e^{-t}\sum_{x=0}^{\infty}\frac{(\theta e^t)^{x+1}}{(x+1)!}\\
&=&e^{\theta e^t}(1-e^{-t})+e^{-t},
\end{eqnarray*}
from which it follows that
\[\e(X)=M^\prime(0)=e^{\theta}-1
\]
and 
\[\v(X)=(2\theta-1)e^{\theta}+1-(e^{\theta}-1)^2=e^{\theta}(2\theta+1-e^{\theta}).\]
Distribution theory inevitably leads into questions of statistical inference for parametric families.  We  might, for example, ask for an estimate of $\theta$ based on a random sample $X_1,\ldots,X_n$ of size $n$.  It turns out that maximum likelihood estimates (``which value of $\theta$ is most likely to have created this data set?") are mathematically intractable and we have to resort to simulations.  A {\it method of moments} (MOM) estimate can easily be found, however.  Here we equate the sample mean $\bar{X}$ and $\e(X)$ to get $\bar{X}=\e(X)=e^{\theta}-1$, so that $\hat{\theta}=\ln(\bar{X}+1)$, where we must have $\bar{X}<e-1$ since $\theta<1$. But, this is not guaranteed, so we ask what is $\p(\bar{X}<e-1)$? We have
\[\p(\bar{X}<e-1)=\p(\frac{\bar{X}-\mu}{\sigma/\sqrt{n}}<\frac{(e-e^{\theta})\sqrt{n}}{e^{\theta/2}\sqrt{2\theta+1-e^{\theta}}}.\] Defining $k_\theta=\frac{e-e^{\theta}}{e^{\theta/2}\sqrt{2\theta+1-e^{\theta}}}$, we see that $k_\theta>0$ for any $\theta\in(0,1)$. The central limit theorem then implies, with $Z$ denoting the standard normal variable,  that $\p(\bar{X}<e-1)\approx \p(Z<\sqrt{n}k_\theta)\to1$ as $n\to\infty$.  Thus, the MOM estimate is reliable for large sample sizes.  In summary, the MOM estimator can be expressed as 
\[ \hat{\theta}_{MOM} = \left\{ 
  \begin{array}{l l}
    \ln(\bar{X}+1) & \quad \text{if $\ln(\bar{X}+1)<1$}\\
    1 & \quad \text{if $\ln(\bar{X}+1)\geq1$}
  \end{array} \right.\]

\noindent 3. THE TELESCOPING GEOMETRIC DISTRIBUTION: The following is a very natural question:  in how many permutations, in which the longest increasing subsequence is of length 2, does the first ascent occur in positions $k, k+1$? Generalizing the numbers that bore his name, Catalan \cite{cata} proved the $k$-fold Catalan convolution formula
\[C_{n,k}:=\sum_{i_1+\ldots+i_k=n}\prod_{r=1}^kC_{i_r-1}=\frac{k}{2n-k}{{2n-k}\choose{n}},\]  and in \cite{cgg} it was proved that there are precisely $C_{n,k}$ permutations on $[n]$ with longest increasing subsequence of size 2 and for which the first ascent occurs at positions $k,k+1$.    Thus, for a randomly chosen 123-avoiding permutation on $[n]$, the distribution of the location of first ascent is given by
\[f(k)=\frac{C_{n,k}}{C_n}=k\frac{(2n-k-1)!(n+1)!}{(2n)!(n-k)!},\enspace k=1,2,\ldots,n,\]
which, for small $k$ and large $n$, may be approximated by
$f(x)=\frac{k}{2^{k+1}}.$
Accordingly, let us define the telescoping geometric-like distribution on ${\bf Z}^+=1,2,\ldots$ by
\[f(x)=\frac{x}{2^{x+1}}=\frac{x+1}{2^x}-\frac{x+2}{2^{x+1}}, x=1,2,\ldots.\]
We see that for a telescoping geometric random variable $X$,
\begin{eqnarray*}\e(X)&=&\sum_{x=1}^{\infty}x\times \frac{x}{2^{x+1}}\\
&=&\frac{1}{8}\sum_{x=1}^{\infty}x(x-1)\lr\frac{1}{2}\rr^{x-2}+\frac{1}{4}\sum_{x=1}^{\infty}x\lr\frac{1}{2}\rr^{x-1}\\
&=&\frac{1}{8}\times\frac{2}{(1-\frac{1}{2})^3}+\frac{1}{4}\times \frac{1}{(1-\frac{1}{2})^2}\\
&=&3.
\end{eqnarray*}
Together with the result from the previous section, we have, roughly speaking, that for a random permutation on a large $[n]$, we expect the first ascent to be at position $e-1\approx 1.718$, whereas this value increases to 3 for a random 123-avoiding permutation. 
Moreover, it makes sense, as with the telescoping Poisson distribution, to define a telescoping $\theta$-analog of the above distribution defined by 
\[
f(x)=\frac{(\theta-1)^2x}{\theta^{x+1}}=\frac{(\theta-1)x+1}{\theta^x}-\frac{(\theta-1)(x+1)+1}{\theta^{x+1}}, x=1,2,\dotsc, \theta>1,
\]
and with m.g.f. given by
\begin{eqnarray*}M(t)=\e(e^{tX})&=&\sum_{x=1}^{\infty}e^{tx}\times \frac{(\theta-1)^2x}{\theta^{x+1}}\\
&=&\frac{(\theta-1)^2e^t}{\theta^2}\sum_{x=1}^{\infty}x\lr\frac{e^t}{\theta}\rr^{x-1}\\
&=&\lr\frac{\theta-1}{\theta-e^t}\rr^2\cdot e^t,~\mbox{if}~\frac{e^t}{\theta}<1,~\mbox{i.e., if}~t<\ln\theta.
\end{eqnarray*}

Launching into inference, if $\theta$ is unknown, we can easily find its maximum likelihood estimate (MLE) by computing the likelihood function $\ell(\theta)$ and then maximizing in the standard fashion: \[
\ell(\theta)=\ell(X_1,\ldots,X_n,\theta)=\frac{(\theta-1)^{2n}\prod_{i=1}^{n}X_i}{\theta^{n+\sum_{i=1}^{n}X_i}},
\]
which yields
\[
\frac{\partial\log\ell(\theta)}{\partial\theta}=\frac{2n}{\theta-1}-\frac{n+\sum_{i=1}^{n}X_i}{\theta}.
\]
Setting $\frac{\partial\log\ell(\theta)}{\partial\theta}=0$, we have $$\frac{2n}{\theta-1}=\frac{n+\sum_{i=1}^{n}X_i}{\theta},$$
which yields the unusual MLE $\hat{\theta}=1+\frac{2}{\bar{X}-1}$.
\par For the method of moments estimation, we verify that  
\[
\e(X)==M'(0)=\frac{2}{\theta-1}+1,\]
and setting $\bar{X}=\frac{2}{\hat{\theta}-1}+1$, we obtain the same MOM estimator $\hat{\theta}=1+\frac{2}{\bar{X}-1}$ as the MLE.

Finally we complete the standard undergraduate Mathematical Statistics agenda by testing simple hypotheses about $\theta$. Consider the null hypothesis $H_0:\theta=\theta_0$ versus the alternative hypothesis $H_1:\theta=\theta_1$. Assume that $\theta_1>\theta_0$.  By the Neyman Pearson theorem, the most powerful test rejects $H_0$ if the likelihood ratio $\frac{\ell(\theta_1)}{\ell(\theta_0)}>k$, where $k>0$ is a constant.  Below, as is customary, we let $k$ be a generic constant whose value changes from line to line:
Since
\[\frac{\ell(\theta_1)}{\ell(\theta_0)}=\lr\frac{\theta_1-1}{\theta_0-1}\rr^{2n}\lr\frac{\theta_0}{\theta_1}\rr^{n+\sum_{i=1}^{n}X_i},
\]
we reject $H_0$ if $$\lr\frac{\theta_1 -1}{\theta_0 -1}\rr^{2n}\lr\frac{\theta_0}{\theta_1}\rr^{n+\sum_{i=1}^{n}X_i}>k,$$ 
which simplifies, on taking logs, to
$$\lr n+\sum_{i=1}^{n}X_i\rr\log\lr\frac{\theta_0}{\theta_1}\rr>k.$$
Since $\theta_1>\theta_0$, we turn the above condition into $$n+\sum_{i=1}^{n}X_i<k,$$
or into the compact
$\bar{X}<k$.  $H_0$ is rejected for small values of $\bar{X}$, where the critical value is determined by the level of significance used.

\medskip

\noindent{\bf Summary.}  Telescoping sums very naturally lead to probability distributions on ${\mathbb Z}^+$.  But are these distributions typically cosmetic and devoid of motivation?  In this paper we give three examples of ``first occurrence" distributions, each defined by telescoping sums, and that each arise from concrete questions about the structure of permutations.

\end{document}